\documentclass[12pt]{article}
\begin{document}
\def\R{\mathbf{R}}
\def\C{\mathbf{C}}
\def\Ima{\mathrm{Im}\, }
\def\Rea{\mathrm{Re}\, }
\title{Stability in the Marcinkiewicz theorem}
\author{Alexandre Eremenko and Alexander Fryntov}
\maketitle

\begin{center}
{\em Dedicated to the memory of I. V. Ostrovskii}
\end{center}

\begin{abstract} Ostrovskii's generalization of
the Marcinkiewicz theorem implies that if an entire
characteristic functions of a probability distribution
satisfies $\log\log M(r,f)=o(r)$
and is zero-free then the distribution is normal.
We show that under the same growth
condition, absence of zeros in a wide vertical strip implies that
the distribution is close to a normal one. This generalizes and simplifies
a recent result of Michelen and Sahasrabudhe.

\end{abstract}

Following Linnik \cite{L}, an entire function
$f$ 
is called  
a {\em ridge function} if $|f(z)|\leq |f(i\Ima z)|,\; z\in \C$. 
This definition
is justified by Probability theory:
characteristic functions of random variables are ridge functions when they are
entire. 
We will apply the same
name to subharmonic functions $u$ in $\C$ satisfying
\begin{equation}
\label{ridge}
u(z)\leq u(i\Ima z),\quad z\in\C.
\end{equation}
Classical theorem of Marcinkiewicz \cite{Ma} says that all ridge entire 
functions of finite order without zeros are of the form $\exp(-az^2+biz+c)$,
where $a>0$, $b$ is real, and $c$ is complex.
This was generalized by Ostrovskii \cite{O}
who proved a conjecture of Linnik that
the condition of finite order can be relaxed to 
$$\log^+\log|f(z)|=o(|z|),\quad z\to\infty.$$
This condition was further relaxed in \cite{K} to
\begin{equation}\label{10}
\liminf_{z\to\infty}\frac{\log^+\log|f(z)|}{|z|}=0.
\end{equation}
Paper \cite{V} contains
a survey of further generalizations of Ostrovskii's result. 

We prove a ``stable version'' of this theorem for entire functions which are
free of zeros in vertical strips:
\vspace{.1in}

\noindent
{\bf Theorem 1.}
{\em If $u$ is a ridge subharmonic function in $\C$ satisfying
\begin{equation}\label{growth}
\liminf_{r\to\infty}\frac{\log\max\{ u(ir),u(-ir)\}}{r}=0,
\end{equation} 
is harmonic in the strip
\begin{equation}
\label{1}
S(\Delta)=\{ z:|\Rea z|<\Delta\}
\end{equation}
and normalized by $u(0)=u_x(0)=u_y(0)=0$ and $u_{yy}(0)=1$,
then 
\begin{equation}
\label{estim}
|u(z)+\Rea(z^2/2)|\leq c_0|z|^3/\Delta,\quad |z|\leq \Delta/3,
\end{equation}
where $c_0$ is an absolute constant.}
\vspace{.1in}

Example $u(z)=\cosh y\cos y-1$ shows that the growth condition (\ref{growth})
is best possible.
A new proof of Linnik's conjecture is obtained
by setting $u=\log|f|$ and $\Delta=\infty$.

As a corollary we obtain a generalization of the recent theorem by
Michelen and Sahasrabudhe
\cite[Thm. 4.1]{M}:
\vspace{.1in}

\noindent
{\bf Theorem 2.} {\em Let $X$ be a random variable with average $\mu$ and
standard
deviation $\sigma$. Suppose that the characteristic function $f_X$
is entire, satisfies (\ref{10}), and is free of zeros in the strip
$\{ z:|\Rea z|<\delta\}$. Then the distribution function $F_{X^*}$
of the random variable $X^*=(X-\mu)/\sigma$ satisfies
$$|F_{X^*}-F_N|_\infty\leq \frac{c_1}{\sigma\delta},$$
where $c_1$ is an absolute constant, and $N$ is the standard normal
distribution with characteristic function $f_N(z)=\exp(-z^2/2)$.} 
\vspace{.1in}

This theorem was proved in \cite{M} under the additional assumption
that $X$ takes values in the set $\{0,1,\ldots,n\}$.
We generalize the result
and propose a shorter proof.
We will use the
\vspace{.1in}

\noindent
{\bf Phragm\'en--Lindel\"of Theorem.} {\em If a subharmonic function $v$
in a strip
$S$
satisfies 
\begin{equation}\label{growth1}
\liminf_{z\to\infty}\frac{\log^+v(z)}{|z|}=0,
\end{equation}
and $v(z)\leq 0,\; z\in\partial S$, then $v(z)\leq 0$ in $S$.}
\vspace{.1in}

\noindent
{\bf Lemma 1.} {\em If a harmonic function in a strip $S(\Delta)$
satisfies (\ref{growth})
and (\ref{ridge}), then for all real $y$, the function
$x\mapsto u(x+iy)$ is decreasing for $x\in[0,\Delta/2]$.}
\vspace{.1in}

{\em Proof.} Let us fix $s\in(0,\Delta/2)$ and let $z\mapsto z^*$
be the reflection with respect to the line $\Rea z=s$, that is
$z^*=2s-\overline{z}$. We define
$u^*(z)=u(z^*)$, and $$v(z)=\max\{ u(z),u^*(z)\},\quad 0<\Rea z<2s.$$

On the lines $\Rea z=0$ and $\Rea z=s$ we have $v(z)\leq u(z)$. For a ridge
function $u$, condition (\ref{growth}) implies (\ref{growth1}) so
$u$ and $v$ satisfy (\ref{growth1}), and  by the
Phragm\'en--Lindel\"of theorem we conclude that  $v(z)\leq u(z)$ in the strip $\{ z:0<\Rea z <s\}$.
On the other hand $v(z)\geq u(z)$ by definition, so 
\begin{equation}\label{a}
v(z)=u(z),\quad 0<\Rea z<s.
\end{equation}

On the lines $\Rea z=s$ and $\Rea z =2s$ we have
$v(z)\leq u^*(z)$, so by a similar application of the
the Phragm\'en--Lindel\"of theorem we conclude that $v(z)\leq u^*(z)$
in the strip $\{ z:s<\Rea z<2s\}$. On the other hand, $v(z)\geq u^*(z)$ by definition, so 
\begin{equation}\label{b}
v(z)=u^*(z),\quad s<\Rea z<2s.
\end{equation}

Since $v(z)$ is subharmonic, we have $v_x(s-0)\leq v_x(s+0)$, and 
in view of (\ref{a}), (\ref{b}) we have
$$v_x(s-0)=u_x(s)\quad\mbox{and}\quad v_x(s+0)=u^*_x(s)=-u_x(s),$$
and so we obtain that $u_x(s)\leq -u_x(s)$ that is $u_x(s)\leq 0$,
which proves the Lemma.
\vspace{.1in}

\noindent
{\bf Lemma 2.} {\em Let  $Q$ be the square,
\begin{equation}
\label{R}
Q=\{ x+iy:0<x<2,\; |y|<1\},
\end{equation}
and
let $P(z,\zeta)$ be the Poisson kernel of $Q$, where $z=x+iy\in Q,$
and $\zeta\in\partial Q$.
Then for $\zeta\in\partial Q\backslash(-i,i)$ we have
$$ P_x(0,\zeta)\geq c_2,$$
where $c_2$ is an absolute constant.}
\vspace{.1in}

\noindent
{\bf Lemma 3.} {\em The family of harmonic functions in a vertical
strip $S(\Delta)$ as in (\ref{1}) satisfying (\ref{growth}), (\ref{ridge}) and
normalized both conditions
$$u(0)=u_y(0)=0,\quad u_{yy}(0)=1,$$
is uniformly bounded from above on every compact set $K\subset S(\Delta/2)$
by a
constant depending only on $K$ and $\Delta$.}

\vspace{.1in}
{\em Proof.} By Lemma 1, harmonic functions $-u_x$ are positive 
in the right half of the strip, and $u_x(0,y)=0$ in view
of (\ref{ridge}). Applying to them the Poisson representation
in rectangles $cQ$ where $Q$ is defined in (\ref{R}) and
using Lemma 2, we obtain that the total measure
in this representation is bounded. So $u_x$ are uniformly bounded
on compacts. We conclude that the analytic functions
$u_x-iu_y$, are uniformly bounded on compacts.
Since $u_x(0)=0$ by the ridge property and $u_y(0)=0$ by assumption,
we conclude that functions $u$ are uniformly bounded on compacts
in $S(\Delta/2)$. 
This proves Lemma 3.
\vspace{.1in}

{\em Proof of Theorem 1.}
We may assume without loss of generality that $\Delta\geq 1.$
Consider the expansion at $0$:
$$u(z)=\Rea \left(-z^2/2+\sum_{n=3}^\infty a_nz^n\right).$$
Let
$$u_\Delta=\Delta^{-2}u(z\Delta)=\Rea\left(-z^2+
\sum_{n=3}^\infty a_n\Delta^{n-2}z^n\right),\quad z\in S(1).$$
By Lemma 3, its coefficients
are uniformly bounded, therefore $|a_n|\leq c_3\Delta^{2-n}$, and
$$\sum_{n=3}^\infty|a_n||z^n|\leq c_3\Delta^{-1}\frac{|z|^3}{1-|z|/\Delta}\leq
c_0|z|^3/\Delta,\quad\mbox{when}\quad|z|\leq\Delta/3.$$
This proves Theorem 1.
\vspace{.1in}

{\em Derivation of Theorem 2 from Theorem 1.} Following \cite{M}
and \cite{R},
we use the Berry--Esseen inequality
\begin{equation}\label{feller}
\sup_{t\in\R}|F_{X^*}(t)-F_Z(t)|\leq\frac{1}{\pi}\int_{-T}^T\left|
\frac{f_{X*}(x)-e^{-x^2/2}}{x}\right|dx+\frac{c}{T},
\end{equation}
where $c$ is an absolute constant.

This estimate can be found in \cite[Ch. XVI, 3, Lemma 2]{Feller} and
in \cite[Lemma 8.2.2]{LO}.

We set $\Delta=\delta\sigma$. The statement of Theorem 2 is meaningful
only when $\Delta$ is large, so we assume that $\Delta>c_0$,
where $c_0$ is the constant in Theorem 1.

We are going to apply Theorem 1 to $u=\log|f_{X^*}|,$ where
$f_{X*}$ is the characteristic function of $X^*$. Since $X^*$ is normalized,
$u$ is normalized as required in Theorem 1. Since by assumption
the characteristic
function $f_X$ has no zeros in the strip $S(\delta)$,
the function $f_{X^*}$ has no zeros in the strip $S(\Delta)$.
Then Theorem 1 implies that
$$f_{X^*}(x)=\exp(-x^2/2+R(x)),\quad\mbox{where}\quad|R(x)|\leq c_0|x|^3/\Delta,\quad|x|<\Delta/2.$$ 
Set $T=\Delta/(4c_0)$ in (\ref{feller}).
To estimate the integral in (\ref{feller}) we break it
into two parts:

Let
$$a:=(\Delta/c_0)^{1/3}\geq 1.$$

When $|x|<a$,
we have $|R(x)|\leq 1$,
so $|e^{R(x)}-1|\leq 2|R(x)|\leq 2c_0x^3/\Delta,$ so
$$\int_{-a}^a\left|\frac{f_{X^*}(x)-e^{-x^2/2}}{x}\right|dx=\frac{2c_0}{\Delta}
\int_{-\infty}^\infty e^{-x^2/2}x^2dx\leq c_5/\Delta.$$

When $|x|\in [a,T]$ we use 
$$f_{X^*}(x)=\exp(-x^2/2+R(x))$$
and
$$x^2(-1/2+|x|c_0/\Delta)\leq x^2(-1/2+1/4)=-x^2/4.$$
So
$$\int_{|x|\in[a,T]}\left|\frac{f_{X*}(x)-e^{-x^2/2}}{x}\right|dx\leq
4\int_a^\infty e^{-x^2/4}dx\leq c_6/\Delta.$$
This completes the proof.

The authors thank F. Nazarov and M. Sodin for useful discussions.

{\em Department of Mathematics,

 Purdue University,

West Lafayette, IN, 47907 USA
\vspace{.1in}

{\em 6198 Townswood ct.,

Mississauga ON, L5N2L4, Canada}

eremenko@math.purdue.edu}
\end{document}